\documentclass{amsart}

\usepackage[utf8]{inputenc}
\usepackage{amsmath,amssymb,amsthm,units, stmaryrd,stackrel,relsize,bm}
\usepackage{tikz}

\usetikzlibrary{arrows,positioning} 
\tikzset{
    >=stealth',
    punkt/.style={
           rectangle,
           rounded corners,
           draw=black, very thick,
           text width=6.5em,
           minimum height=2em,
           text centered},
    pil/.style={
           ->,
           thick,
           shorten <=2pt,
           shorten >=2pt,}
}
\usetikzlibrary{positioning}

\newtheorem{theorem}{Theorem}
\newtheorem{conjecture}[theorem]{Conjecture}
\newtheorem{lemma}[theorem]{Lemma}

\newtheorem{claim}[theorem]{Claim}
\newtheorem{corollary}[theorem]{Corollary}
\newtheorem{proposition}[theorem]{Proposition}

\theoremstyle{definition}
\newtheorem{definition}[theorem]{Definition}

\theoremstyle{remark}
\newtheorem{remark}[theorem]{Remark}

\usepackage[
    colorlinks=true, 
    citecolor=red,
    linkcolor=blue,
    urlcolor=blue]{hyperref}
\newcommand\ad{{\mathsf{AD}}}

\newcommand\hod{{\mathsf{HOD}}}
\newcommand\zfc{{\mathsf{ZFC}}}
\newcommand\zf{{\mathsf{ZF}}}
\newcommand\kp{{\mathsf{KP}}}
\newcommand\dc{{\mathsf{DC}}}

\newcommand\Ord{{\mathsf{Ord}}}

\newcommand\coll{\text{Coll}}

\newcommand{\SIGMA}{\boldsymbol\Sigma}
\newcommand{\DELTA}{\boldsymbol\Delta}

\newcommand{\ind}{\textsf{Ind}}
\newcommand{\dGame}{\breve{\Game}}

\makeatletter
\def\Ddots{\mathinner{\mkern1mu\raise\p@
\vbox{\kern7\p@\hbox{.}}\mkern2mu
\raise4\p@\hbox{.}\mkern2mu\raise7\p@\hbox{.}\mkern1mu}}
\makeatother

\begin{document}
\title{Determined Admissible Sets}
\author{J. P. Aguilera}
\address{Department of Mathematics, University of Ghent. Krijgslaan 281-S8, B9000 Ghent, Belgium.}
\address{Institute of Discrete Mathematics and Geometry, Vienna University of Technology. Wiedner Hauptstra{\ss}e 8--10, 1040 Vienna, Austria.}
\email{aguilera@logic.at}

\begin{abstract}
It is shown, from hypotheses in the region of $\omega^2$ Woodin cardinals, that there is a transitive model of Kripke-Platek set theory containing $\mathbb{R}$ in which all games on $\mathbb{R}$ are determined. 
\end{abstract}
\date{\today}
\clearpage
\maketitle

\setcounter{tocdepth}{1}
\section{Introduction}
Given sets $X$ and $A\subset X^\mathbb{N}$, the Gale-Stewart game of length $\omega$ on $A$ is defined as follows: two players, I and II, alternate turns playing elements of $X$ infinitely many times, producing a sequence $x \in X^\mathbb{N}$. Player I wins if $x \in A$; otherwise, Player II wins. A \emph{strategy} is a function $\sigma:\bigcup_{n\in\mathbb{N}}X^{n} \to X$. We say that $\sigma$ is a \emph{winning} strategy for Player I if $x \in A$ whenever $x(2n) = \sigma(x(0),x(1),\hdots, x(2n-1))$ for all $n\in\mathbb{N}$, i.e., if Player I wins every run of the game obtained by obeying $\sigma$. The notion of  winning strategy for Player II is defined analogously.
$A$ is \emph{determined} if one of the players has a winning strategy for this game. 

The Axiom of Determinacy ($\ad$) states that all Gale-Stewart games of length $\omega$ on (subsets of) $\mathbb{N}$ are determined. Similarly, the Axiom of Real Determinacy ($\ad_\mathbb{R}$) states that all Gale-Stewart games of length $\omega$ on $\mathbb{R}$ are determined. These axioms are very powerful, and have a great deal of consequences on the structure of the real numbers. For instance, they imply that every subset of $\mathbb{R}$ is Lebesgue-measurable, has the property of Baire, and is either countable or contains a perfect subset (see e.g., Kanamori \cite{Ka08}); it follows that both axioms are inconsistent with the Axiom of Choice. By work of Woodin, they are consistent with Zermelo-Fraenkel set theory ($\zf$) without the Axiom of Choice (see Woodin \cite{Wo88} for a proof of the consistency of $\zf + \ad$), although this fact is not itself provable in $\zfc$. To do this, one needs to make use of \emph{large cardinal axioms}, strengthenings of the axiom of infinity. A theorem of Woodin (see Koellner-Woodin \cite{KW10} and Steel \cite{Sta}) states that $\zf +\ad$ is consistent if, and only if, $\zfc$ is consistent with the existence of infinitely many Woodin cardinals. Another theorem of Woodin (unpublished, but see the remark following Steel \cite[Theorem 2.11]{St08d}) shows that $\zf + \ad_\mathbb{R}$ is consistent if, and only if, $\zfc$ is consistent with the existence of a cardinal $\lambda$ which is a limit of Woodin cardinals and of cardinals which are ${<}\lambda$-strong. In particular, it is stronger than $\zfc$ with a proper class of Woodin cardinals, and much stronger than $\zf + \ad$.

Kripke-Platek set theory ($\kp$) is a weak set theory studied for its recursion-theoretic properties. Many results about recursion on $\mathbb{N}$ generalize to arbitrary transitive sets which satisfy all axioms of $\kp$; these are called \emph{admissible sets}. A natural question is that of the strength of the theories $\kp+\ad$ and $\kp+\ad_\mathbb{R}$ (here, $\ad$ and $\ad_\mathbb{R}$ are assumed to include a clause asserting the existence of $\mathcal{P}(\mathbb{N})$, and thus of the real line -- this is not provable in $\kp$). Clearly, $\zf+\ad$ proves the consistency of $\kp+\ad$, although the latter theory is not significantly weaker than the former; e.g., $\kp+\ad$ proves the consistency of $\zfc$ with all finite amounts of Woodin cardinals, and much more (see the remark on p. \pageref{theremark}). One may be led to conjecture  that $\kp+\ad_\mathbb{R}$ has similar consistency strength to that of $\zf+\ad_\mathbb{R}$. Here, we observe that this is not the case at all:\\

\paragraph{\textbf{Theorem}}
\textit{
Suppose there are $\omega^2$ Woodin cardinals and a measurable cardinal above. Then, there is a transitive model of $\kp + \ad_\mathbb{R}$ containing all reals.}\\

The conclusion of the theorem is in fact proved under a weaker hypothesis in Section \ref{Sectmainopencubed} (Theorem \ref{mainopencubed}). We conjecture that this hypothesis  is optimal. 
The theorem as stated is proved in Section \ref{Sectweaker}, after which we also sketch an argument for eliminating the measurable cardinal.

The discrepancy between Woodin's theorem on $\ad_\mathbb{R}$ and Theorem \ref{mainopencubed} can be explained by recalling Woodin's (unpublished) theorem whereby $\ad_\mathbb{R}$ is equivalent (over $\zf$ + $\ad$ + $\dc$) to the assertion that all sets of reals can be uniformized. The model $M$ considered in the proof below will satisfy that all sets can be uniformized.
However, even if one assumes that $M$ can be extended to a transitive model of $\zf+\ad$, the least such extension $M'$ will have many sets that cannot be uniformized \emph{within $M'$}.

The main tool for the proof is the theory of inductive definitions with quantifiers on $\mathbb{R}$ and Spector classes, the basics of which we recall in the following section. The reader may consult Aczel \cite{Acz75},
Kechris \cite{Ke16}, or Moschovakis \cite{Mo74} for further background, as well as some of the arguments. As is customary, we identify the real numbers $\mathbb{R}$ and the Baire space $\mathbb{N}^\mathbb{N}$ in what follows. This shall cause no confusion, since the usual real line does not play a role in the argument. One of the main advantages of this identification is that, with it, $\mathbb{R}$ is recursively homeomorphic to the spaces $\mathbb{R}^n ({=} \mathbb{N}^{\mathbb{N}\times n})$ and $\mathbb{R}^\mathbb{N} ({=}\mathbb{N}^{\mathbb{N}\times\mathbb{N}})$, and this allows us to further identify $\mathbb{R}$ with its finite and $\mathbb{N}$-indexed products and, in particular, think of the winning set of a game on $\mathbb{R}$ as a subset of $\mathbb{R}$. We will do such frequently and without mention in order to simplify notation, unless we think it may cause confusion, in which case we state the space explicitly. We shall also refer to sets of pairs of reals as subsets of $\mathbb{R}$, etc.

It might be worth noting that the existence of the real line is equivalent over $\kp$ to the existence of $\mathbb{N}^\mathbb{N}$ and to that of $\mathcal{P}(\mathbb{R})$. Moreover, since the Baire space $\mathbb{N}^\mathbb{N}$ has the same cardinality as the real line, provably in $\kp$, the axioms $\ad_\mathbb{R}$ and $\ad_{\mathbb{N}^\mathbb{N}}$ are equivalent. Hence, these conventions do not affect the result proved.\\

\paragraph{\textbf{Acknowledgements}} The author would like to thank both reviewers for their many comments and suggestions, which certainly improved the article. This work was partially supported by FWF grants P-31063 and P-31955.

\section{Preliminaries}\label{SectPrel}

A (non-trivial) \emph{quantifier on $\mathbb{R}$} is a nonempty collection of subsets of $\mathbb{R}$ which is closed under supersets and not equal to $\mathcal{P}(\mathbb{R})$.
If $Q$ is a quantifier on $\mathbb{R}$, we write $Qx\, X(x)$ for $X \in Q$ and define its dual by
$\breve{Q}x\, X(x)\leftrightarrow \lnot Qx\, \lnot X(x)$.
Examples of quantifiers in which we will be interested are
$$\exists^\mathbb{R}  = \{A\subset\mathbb{R}: A \text{ is nonempty}\},$$
and
$$\Game^\mathbb{R}  = \{A\subset\mathbb{R}: \text{Player I has a winning strategy the game on $\mathbb{R}$ with payoff $A$}\}.$$
We may write $\exists$ for $\exists^\mathbb{R}$ if it is clear from the context that the quantification is carried out over reals.

Fix a quantifier $Q$ on $\mathbb{R}$. Although our main case of interest is $Q = \Game^\mathbb{R}$, we will need other examples, so we shall work in the abstract in this section.
Since $Q$ is closed under supersets, $X\in Q$ is equivalent to $\exists Y \in Q\, \forall x\in Y\, X(x)$. Motivated by this, one can generalize formulae of the form $Qx\, X(x)$ to allow for infinite strings of quantifiers: 
\[Q_0 x_0\, Q_1 x_1 \, \hdots X(x_0,x_1,\hdots).\]
The displayed formula is \emph{defined} to hold if Player I has a winning strategy in the following game: during turn $k$, Player I chooses $X_k \in Q_k$ and Player II responds with $x_k \in X_k$. After infinitely many rounds have been played, Player I wins if, and only if,
\[(x_0, x_1, \hdots) \in X.\]
The reader might want to verify that the expected equality
\[\Game^\mathbb{R} x\, X(x) \leftrightarrow \exists x_0\, \forall x_1\, \exists x_2\,\hdots X(x_0,x_1,x_2,\hdots)\]
holds under these definitions (recall that we are identifying the spaces $\mathbb{R}$ and $\mathbb{R}^\mathbb{N}$).\,

For our purposes, an \emph{operator} is a function
\[\phi:\mathcal{P}(\mathbb{R})\to\mathcal{P}(\mathbb{R}).\]
If $\phi$ is an operator, then it can be iterated by setting
\begin{align*}
\phi^0 
&= \varnothing\\
\phi^{<\alpha} 
&= \bigcup_{\beta<\alpha}\phi^\beta\\
\phi^\alpha &= \phi(\phi^{<\alpha})\\
\phi^\infty &= \bigcup_{\alpha \in \Ord}\phi^\alpha.
\end{align*}
Consider the language $L(Q)$ of second-order arithmetic augmented with symbols for $Q$ and $\breve{Q}$ as quantifiers, constants for each element of $\mathbb{R}$, and a predicate symbol $\dot X$ for a set of reals. 
 
Given a formula $\psi$ in $L(Q)$, we say that $\dot X$ appears \emph{only positively} in $\psi$ if $\psi$ is built up by formulae not involving $\dot X$ and atomic formulae of the form $\dot X(x)$ using the connectives $\wedge$ and $\vee$ and quantifiers. 
We are interested in operators that are definable by a closed (i.e., variable-free) formula in $L(Q)$ in which $\dot X$ appears only positively, i.e., operators of the form 
\[X \mapsto \{x\in\mathbb{R}: \psi(x) \text{ holds when $\dot X$ is interpreted as $X$}\},\]
for some closed $L(Q)$-formula $\psi$ in which $\dot X$ appears only positively. If $\phi$ is such an operator, we say that $\phi^\infty$ is the subset of $\mathbb{R}$ \emph{inductively defined} by $\phi$. If so, we denote by $|\phi|$ the \emph{closure ordinal} of $\phi$: the least ordinal $\alpha$ such that $\phi^{<\alpha} = \phi^\infty$. 

We write $\ind(Q)$ for the pointclass of all sets of the form
\[\big\{x : (x, a) \in \phi^\infty\big\}\]
for some $a\in\mathbb{R}$ and some
operator $\phi$ which is definable by a closed formula in $L(Q)$ in which $\dot X$ appears only positively.
We also say:
$$\text{$A$ is $Q$-hyperprojective} \longleftrightarrow A \in \ind(Q) \wedge \mathbb{R}\setminus A \in \ind(Q).$$
We define $|\ind(Q)|$ to be the supremum of $|\phi|$, where $\phi$ is an operator as above. 
$|\ind(Q)|$ is also the supremum of lengths of prewellorderings coded by $Q$-hyperprojective sets. Given $A \in \ind(Q)$, as witnessed by an operator $\phi$ and a parameter $a\in\mathbb{R}$, Moschovakis' \emph{Stage Comparison} Theorem 
states that there is a norm in $\ind(Q)$ on $A$ which is isomorphic to the function
\[x \mapsto \min\{\alpha : (x,a)\in \phi^{\alpha}\}.\]
(It might not be the displayed function itself, as it might not be surjective for some parameters $a$.)
The \emph{Boundedness} Theorem 
states that if $A$, $\phi$ and $a$ are as above, then $A$ is $Q$-hyperprojective if, and only if, $|\phi| < |\ind(Q)|$.

It is customary to employ quantifiers also as operators from $\mathcal{P}(\mathbb{R})$ to $\mathcal{P}(\mathbb{R})$. With this (abuse of) notation, $Q\, A$ denotes the collection of all $y\in\mathbb{R}$ such that $\{x : (x,y) \in A\}$ belongs to $Q$, in the sense defined above.\\

A \emph{Spector class} on $\mathbb{R}$ is an $\mathbb{R}$-parametrized collection of subsets of $\mathbb{R}$ closed under $\wedge,\vee,\exists^\mathbb{R},\forall^\mathbb{R}$, containing all projective sets and having the prewellordering property. We will not make direct use of this definition, but we will need some facts about Spector classes. By a theorem of Aczel \cite{Acz75}, $\ind(Q)$ is the smallest Spector class on $\mathbb{R}$ closed under $Q$ and $\breve{Q}$, in the sense that
\[A \in \ind(Q) \longrightarrow Q\,A \in \ind(Q) \wedge \breve{Q}\,A \in \ind(Q).\]
Suppose that $A$ is $Q$-hyperprojective, so that $A \in \ind(Q)$ and $\mathbb{R}\setminus A \in\ind(Q)$. By Aczel's theorem, $Q\, A\in \ind(Q)$ and $\breve{Q}\,(\mathbb{R}\setminus A) \in \ind(Q)$. Since
$\breve{Q}(\mathbb{R}\setminus A) = \mathbb{R}\setminus Q\, A$, we see that $Q\, A$ is also $Q$-hyperprojective, so that the $Q$-hyperprojective sets are closed under $Q$ and $\breve{Q}$. Since existential and universal quantifiers over $\mathbb{R}$ come for free with the language $L(Q)$, the $Q$-hyperprojective sets are also closed under $\exists^\mathbb{R}$ and $\forall^\mathbb{R}$. They are also closed under continuous preimages, and therefore under countable unions and intersections.

Another theorem of Aczel \cite{Acz75} gives an alternate characterization of $\ind(Q)$ 
\[\text{$\ind(Q) = \{Q^+\, A: A$ is projective$\}$},\]
where $Q^+$ denotes the \emph{next quantifier} of $Q$, given by 
$$Q^+u\, R(u) \leftrightarrow Qx_0\,\breve{Q}x_1\,\exists x_2\, \forall x_3\, Qx_4\, \hdots\, \exists n\, R(\langle x_0, x_1,\hdots, x_n \rangle),$$
where the infinite strings of quantifiers is interpreted as above.

\section{A model of $\kp + \ad_\mathbb{R}$} \label{Sectmainopencubed}
The same way one defines games of length $\omega$ on $\mathbb{N}$ or $\mathbb{R}$, one can define games of transfinite length. For definiteness, we assume that limit moves are made by Player I, though this choice will be immaterial for the proof. 

\begin{theorem} \label{mainopencubed}
Suppose that open games of length $\omega^3$ on $\mathbb{N}$ are determined. Then, there is a transitive model of $\kp + \ad_\mathbb{R}$ containing all reals.
\end{theorem}

The first step to proving Theorem \ref{mainopencubed} is to consider the quantifier
\begin{align*}
\Game_{\omega^2}^\mathbb{R}   
&= \big\{A\subset\mathbb{R}: \text{Player I has a winning strategy the game }\\
& \quad \quad \text{ of length $\omega^2$ on $\mathbb{R}$ with payoff $A$}\big\}.
\end{align*}
Using the quantifiers-as-operators notation, we have
\[\Game_{\omega^2}^\mathbb{R} A = \Big\{y\in\mathbb{R} : \{x: (x,y) \in A\} \in \Game_{\omega^2}^\mathbb{R}\Big\}.\]
We write
\begin{equation*}\label{eqOpen}
\Game^\mathbb{R}_{\omega^2}\SIGMA^0_1 = \big\{\Game^\mathbb{R}_{\omega^2}A: A \text{ is open in the usual topology}\big\}
\end{equation*}
We state two lemmata.
\begin{lemma}\label{LemmaOpenCubed1}
$\ind(\Game^\mathbb{R}) \subset \Game^\mathbb{R}_{\omega^2}\SIGMA^0_1$.
\end{lemma}

\begin{lemma} \label{LemmaOpenCubed3}
Suppose that $\Game^\mathbb{R}$-hyperprojective games of length $\omega^2$ on $\mathbb{N}$ are determined. Then every $\Game^\mathbb{R}$-hyperprojective game of length $\omega$ on $\mathbb{R}$ has a $\Game^\mathbb{R}$-hyperprojective winning strategy.
\end{lemma}

Before proving Lemma \ref{LemmaOpenCubed1} and Lemma \ref{LemmaOpenCubed3}, let use show how Theorem \ref{mainopencubed} is a consequence of them:

\begin{lemma}\label{LemmaOpenomega3}
Suppose that open games of length $\omega^3$ on $\mathbb{N}$ are determined. Then, games of length $\omega^2$ on $\mathbb{N}$ with payoff in $\Game^\mathbb{R}_{\omega^2}\SIGMA^0_1$ are determined.
\end{lemma}
\proof
Consider a game of length $\omega^2$ on $\mathbb{N}$ with payoff $A \in \Game^\mathbb{R}_{\omega^2}\SIGMA^0_1$. By the definition of $\Game^\mathbb{R}_{\omega^2}\SIGMA^0_1$, there is an open $B\subset\mathbb{R}$ such that for all $x\in\mathbb{N}^{\omega^2}$, $x\in A$ if, and only if, Player I has a winning strategy in the game of length $\omega^2$ on $\mathbb{R}$ with payoff (the subset of $\mathbb{R}^{\omega^2}$ identified with)
\[\{y\in\mathbb{R}: (x,y)\in B\}.\]
We consider the following game of length $\omega^3$:  
\begin{enumerate}
\item Players I and II begin alternating $\omega^2$-many rounds playing natural numbers to produce a sequence $x \in\mathbb{N}^{\omega^2}$.
\item Afterwards, Players I and II alternate $\omega^2$-many blocks of $\omega$-many turns each. In each block, the corresponding player plays all the digits of a real number $y_i$, for $i < \omega^2$.
\item At the end of the game, Player I wins if $(x, \langle y_i : i\in\omega^2\rangle) \in B$.
\end{enumerate}
The length of this game is $\omega^2 + \omega^2 \cdot \omega = \omega^3$, and its payoff is open, so it is determined. Given a winning strategy $\Sigma$ (for either player), it is easy to verify that the restriction of $\Sigma$ to the first $\omega^2$-many moves is a winning strategy for the game of length $\omega^2$ with payoff $A$.
\endproof

\proof[Proof of Theorem \ref{mainopencubed}]
By the first of the two theorems of Aczel mentioned at the end of  the previous section, $\ind(\Game^\mathbb{R})$ is a Spector class on $\mathbb{R}$. The Companion Theorem of Moschovakis \cite[Theorem 9E.1]{Mo74} then implies that there is an admissible set $M$ (the \emph{companion} of $\ind(\Game^\mathbb{R})$) whose sets of reals are precisely the $\Game^\mathbb{R}$-hyperprojective sets.

Suppose that open games of length $\omega^3$ are determined.
By Lemma \ref{LemmaOpenomega3}, games of length $\omega^2$ on $\mathbb{N}$ with payoff in $\Game^\mathbb{R}_{\omega^2}\SIGMA^0_1$ are also determined. 
By Lemma \ref{LemmaOpenCubed1}, $\Game^\mathbb{R}$-hyperprojective games of length $\omega^2$ on $\mathbb{N}$ are determined. By Lemma \ref{LemmaOpenCubed3}, every $\Game^\mathbb{R}$-hyperprojective game of length $\omega$ on $\mathbb{R}$ has a $\Game^\mathbb{R}$-hyperprojective winning strategy, so $M$ satisfies $\ad_\mathbb{R}$. 
\endproof

\begin{remark}\label{theremark}
The converse of Lemma \ref{LemmaOpenCubed1} is true and can be proved using Aczel's characterization of $\ind(Q)$. The converse of Lemma \ref{LemmaOpenomega3} is also true; it can be proved by the argument of \cite[Section 3]{Ag18b} (having Player II be the closed player).\qed
\end{remark}
\begin{remark}
Woodin has shown that, over $\zf + \dc$, $\ad_\mathbb{R}$ implies that every game of fixed countable length with moves in $\mathbb{N}$ is determined. As the interested reader will verify using the model $M$ above, $\zf$ cannot be replaced by $\kp$ in the statement of Woodin's theorem.\qed
\end{remark}

\begin{remark} \label{remarkad}
A similar argument to the one above shows that if open games of length $\omega^2$ are determined, then there is a transitive model of $\kp + \ad$ containing all reals. In this case, the converse can be proved by using the Kechris-Woodin determinacy transfer theorem \cite{KW83}. \qed
\end{remark}

In order to complete the proof of the theorem, it remains to prove Lemma \ref{LemmaOpenCubed1} and Lemma \ref{LemmaOpenCubed3}.
Given two strategies $\sigma$ and $\tau$, we denote by $\sigma*\tau$ the result of facing them off against each other. In case these are strategies for games on $\mathbb{R}$, we might abuse notation and identify $\sigma*\tau$ with the real coding the play.
Given a strategy $\sigma$ for a player and a sequence of moves $x$ for the opponent, we denote by $\sigma(x)$ the response to $x$ given by $\sigma$.
\proof[Proof of Lemma \ref{LemmaOpenCubed1}]
We need to show that $\ind(\Game^\mathbb{R}) \subset \Game^\mathbb{R}_{\omega^2}\SIGMA^0_1$.
We emphasize that the lemma is stated with no determinacy assumptions.

Suppose $A$ belongs to $\ind(\Game^\mathbb{R})$. By the second of the two theorems of Aczel mentioned at the end of the previous section, there is an analytical\footnote{i.e., parameter-free first-order formula with all quantifiers bounded by $\mathbb{R}$.} $\phi$ and some sequence $\vec a$ of parameters in $\mathbb{R}$ such that
$$x \in A \leftrightarrow (\Game^\mathbb{R})^+u\, \phi(u, x, \vec a).$$
By definition of $(\Game^\mathbb{R})^+$,
\begin{equation}\label{eqInd1}
x\in A \longleftrightarrow \Game^\mathbb{R}x_0\,\dGame^\mathbb{R}x_1\, \exists x_2\, \forall x_3\, \Game^\mathbb{R} x_4\, \hdots\, \exists n\, \phi(\langle x_0,\hdots, x_n\rangle, x, \vec a).
\end{equation}
Observe that each of the quantifiers $\Game^\mathbb{R}$ and $\dGame^\mathbb{R}$ subsumes each of $\exists$ and $\forall$, so, by replacing $\phi$ with a different formula if necessary, equation \eqref{eqInd1} can be rewritten in the form 
\begin{equation}
x\in A \longleftrightarrow \Game^\mathbb{R}x_0\,\dGame^\mathbb{R}x_1\, \Game^\mathbb{R} x_2\, \dGame^\mathbb{R} x_3\, \Game^\mathbb{R} x_4\, \hdots\, \exists n\, \phi(\langle x_0,\hdots, x_n\rangle, x, \vec a),
\end{equation}
which is slightly more uniform.
By definition, $x\in A$ if, and only if, Player I has a winning strategy in the following game $\mathcal{G}(x)$: During turn $k$, Player I chooses $X_k \in \Game^\mathbb{R}$ if $k$ is even, or $X_k \in \dGame^\mathbb{R}$ if $k$ is odd, and Player II responds with $x_k \in X_k$. After infinitely many rounds have been played, Player I wins if, and only if, 
\[\exists n\, \phi(\langle x_0,\hdots x_n\rangle,x,\vec a).\]
In order to show that $A \in \Game^\mathbb{R}_{\omega^2}\SIGMA^0_1$, we will define two auxiliary games $G(x)$ and $H(x)$. 
It will be easy to show that $\mathcal{G}(x)$ is equivalent to $G(x)$ and that $G(x)$ is equivalent to $H(x)$, and $H(x)$ will be an open game (in the topology induced via a homeomorphism between $\mathbb{R}$ and $\mathbb{R}^{\omega^2}$) of length $\omega^2$ with moves on $\mathbb{R}$.  Moreover, the games $G(x)$ and $H(x)$ will be defined uniformly for all $x \in\mathbb{R}$. 

The game $G(x)$ is defined as follows:
\begin{enumerate}
\item During a turn $k$ of one of the forms $4n$ or $4n+3$, Player I plays a strategy $\sigma_{k}$ for a game of length $\omega$ on $\mathbb{R}$;
\item During a turn $k$ of one of the forms $4n+1$ or $4n+2$, Player II plays a strategy $\sigma_k$ for a game of length $\omega$ on $\mathbb{R}$.
\end{enumerate}
Letting $x_{2k}$ be (the real coding) $\sigma_{4k}*\sigma_{4k+1}$ and $x_{2k+1}$ be (the real coding) $\sigma_{4k+3}*\sigma_{4k+2}$, Player I wins if, and only if there is $n$ such that $\phi(\langle x_0,\hdots, x_n\rangle, x, \vec a)$ holds.

\begin{claim}
Player I has a winning strategy in $\mathcal{G}(x)$ if, and only if, she has one in $G(x)$.
\end{claim}
\proof

Suppose Player I has a winning strategy $\Sigma$ for $\mathcal{G}(x)$. We use $\Sigma$ to describe a winning strategy for $G(x)$.
The strategy consists of playing a run $p$ of the game $G(x)$ while at the same time appealing  to an imaginary run $q$ of the game $\mathcal{G}(x)$ as follows:

At the beginning of the game, $\Sigma$ dictates that Player I should play some set $X_0\in \Game^\mathbb{R}$. By definition, the Gale-Stewart game on $X_0$ with moves in $\mathbb{R}$ is won by Player I, as witnessed by some strategy $\sigma_0$. Choose $\sigma_0$ as the first play of $p$. In $G(x)$, Player II responds with a strategy $\sigma_1$, giving rise to a real
$x_0 = \sigma_0*\sigma_1.$
By the choice of $\sigma_0$, we must have $x_0 \in X_0$, thus making this a valid response to $X_0$ by Player II in $q$.
$\Sigma$ dictates that Player I should play some set $X_1 \in \dGame^\mathbb{R}$. By definition, the Gale-Stewart game on $X_1$ with moves in $\mathbb{R}$ is not won by Player II. Thus, imagine that Player II continues the play $p$ in $G$ with a strategy $\sigma_2$. Since $\sigma_2$ cannot be a winning strategy for the game $X_1$, there is some strategy $\sigma_3$ such that, letting
$x_1 = \sigma_3*\sigma_2,$
we have $x_1 \in X_1$, thus making $x_1$ a valid response for Player II in $q$. We continue simultaneously defining $p$ and $q$ following the diagram below:\\

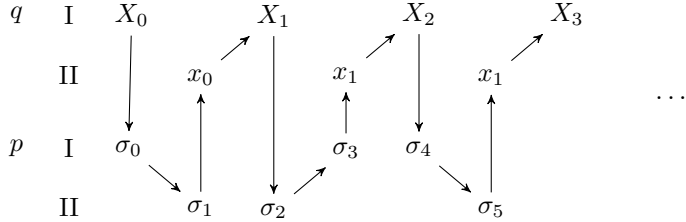
\begin{tikzpicture}[node distance=0.3cm, auto,]

\node (q) {$q$};
\node[below=1.3cm of q] (p) {$p$};
\node[right=of q] (qI) {I};
\node[below=of qI] (pI) {II};
\node[right=of p] (pI) {I};
\node[below=of pI] (pII) {II};

\node[right=of qI] (X0) {$X_0$};

\node[right=of pI] (s0) {$\sigma_0$}
edge[<-] (X0);
\node[below right=0.5cm of s0] (s1) {$\sigma_1$}
edge[<-] (s0);
\node[above=1.3cm of s1] (x0) {$x_0$}
edge[<-] (s1);
\node[above right=0.45cm of x0] (X1) {$X_1$}
edge[<-] (x0);
\node[below=2.1cm of X1] (s2) {$\sigma_2$}
edge[<-] (X1);
\node[above right=0.5cm of s2] (s3) {$\sigma_3$}
edge[<-] (s2);
\node[above=0.55cm of s3] (x1) {$x_1$}
edge[<-] (s3);
\node[above right=0.45cm of x1] (X2) {$X_2$}
edge[<-] (x1);
\node[below=1.3cm of X2] (s4) {$\sigma_4$}
edge[<-] (X2);
\node[below right=0.5cm of s4] (s5) {$\sigma_5$}
edge[<-] (s4);
\node[above=1.3cm of s5] (x2) {$x_1$}
edge[<-] (s5);
\node[above right=0.5cm of x2] (X2) {$X_3$}
edge[<-] (x2);
\node[below right=1cm of X2] (dots) {$\dots$};
\end{tikzpicture}\\

The argument above shows that this results in a play $q$ consistent with $\Sigma$ which is therefore won by Player I. However, $q$ being won by Player I depends only on $\phi$ and $\langle x_0,x_1,\hdots \rangle$, as does $p$ being won by Player I, so $p$ is also won by Player I, as desired.

The converse -- that if Player I has a winning strategy in $G(x)$, then she has one in $\mathcal{G}(x)$ -- is proved similarly. The diagram for this argument is as follows, though in this case some of the arrows do not have a meaning as straightforward as they did before.

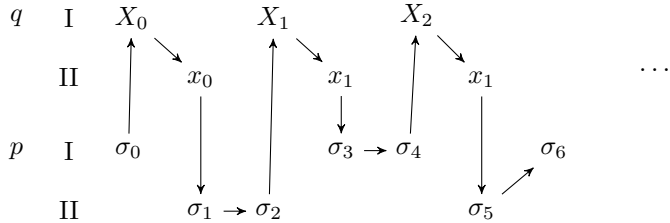
\begin{tikzpicture}[node distance=0.3cm, auto,]
\node (q) {$q$};
\node[below=1.3cm of q] (p) {$p$};
\node[right=of q] (qI) {I};
\node[below=of qI] (pI) {II};
\node[right=of p] (pI) {I};
\node[below=of pI] (pII) {II};

\node[right=of qI] (X0) {$X_0$};

\node[right=of pI] (s0) {$\sigma_0$}
edge[->] (X0);
\node[below right=0.5cm of s0] (s1) {$\sigma_1$};
\node[above=1.3cm of s1] (x0) {$x_0$}
edge[<-] (X0)
edge[->] (s1);
\node[above right=0.45cm of x0] (X1) {$X_1$};
\node[right=of s1] (s2) {$\sigma_2$}
edge[<-] (s1)
edge[->] (X1);
\node[above right=0.5cm of s2] (s3) {$\sigma_3$};
\node[above=0.5cm of s3] (x1) {$x_1$}
edge[<-] (X1)
edge[->] (s3);
\node[above right=0.5cm of x1] (X2) {$X_2$};
\node[right=of s3] (s4) {$\sigma_4$}
edge[<-] (s3)
edge[->] (X2);
\node[below right=0.5cm of s4] (s5) {$\sigma_5$};
\node[above=1.3cm of s5] (x2) {$x_1$}
edge[<-] (X2)
edge[->] (s5);
\node[above right=0.5cm of s5] (X2) {$\sigma_6$}
edge[<-] (s5);
\node[above right=1cm of X2] (dots) {$\dots$};
\end{tikzpicture}\\
Let us explain the diagram.
Fix a winning strategy $\Sigma$ for Player I in $G(x)$; we use it, together with the input from Player II, to construct plays $p$ and $q$, where $p$ is consistent with $\Sigma$, and $q$ is a winning play in $\mathcal{G}(x)$.
To begin, $\sigma_0$ is provided by $\Sigma$. This is used to generate 
\[X_0 = \{x\in\mathbb{R}: x = \sigma_0*\sigma\text{ for some strategy $\sigma$}\},\]
which belongs to $\Game^\mathbb{R}$, as witnessed by the strategy $\sigma_0$.
Afterwards, $x_0$ is provided by Player II in $q$ and $\sigma_1$ is a strategy selected so that 
\[x_0 = \sigma_0 *\sigma_1;\]
this must exist by the choice of $X_0$.
So far, $\sigma_1$ is the only move made by Player II in $p$. If Player II proceeds by playing some $\sigma_2$, then the strategy $\Sigma$ provides an answer which following our earlier convention is denoted by $\Sigma(\sigma_1,\sigma_2)$.
To begin the second turn in $q$, let
\[X_1 = \{x\in\mathbb{R}: x = \Sigma(\sigma_1,\sigma)* \sigma \text{ for some strategy $\sigma$}\},\]
i.e., $X_1$ is the collection of all reals $\langle y_0, y_1,\hdots\rangle$ coding plays that arise from facing off some strategy $\sigma$ for Player II against the strategy for Player I that would be provided by $\Sigma$ if Player II were to play $\sigma$ in $p$. Observe that Player II cannot have a winning strategy for the Gale-Stewart game on reals with payoff $X_1$, for given any strategy $\sigma$ for Player II, $\Sigma(\sigma_1,\sigma)$ is a strategy for Player I that defeats $\sigma$ in this game; thus, $X_1 \in \dGame^\mathbb{R}$, making it a legal move in $q$. Let $x_1$ be Player II's response to $X_1$. Since $x_1 \in X_1$, there is some strategy $\sigma_2$ such that
\[x_1 = \sigma_3*\sigma_2,\]
where $\sigma_3 = \Sigma(\sigma_1,\sigma_2)$. The remainder of the argument is as before. This proves the claim.
\endproof

We have shown that $x\in A$ if, and only if, Player I has a winning strategy in $G(x)$. This game is given in terms of the analytical formula $\phi$ and some parameter $\vec a$.
For simplicity, let us assume that $\phi$ is $\Sigma^1_1$, so that
$$\phi(u,x,\vec a) \leftrightarrow \exists y\in\mathbb{R}\,\forall m\, R(u\upharpoonright m,x\upharpoonright m,\vec a\upharpoonright m, y\upharpoonright m),$$
for some recursive relation $R$.

Let us define a game $H(x)$ of length $\omega^2$ with moves on $\mathbb{R}$ and open payoff which is equivalent to $G(x)$. The game consists of infinitely many blocks, each of length $\omega$, and has two stages.
\begin{enumerate}
\item During turns $(\omega\cdot n, \omega\cdot(n+1))$, if the game is in the first stage, then players I and II alternate turns playing real numbers $x^n_0, x^n_1,\hdots$, etc. If $n$ is even, then Player I begins; if $n$ is odd, then Player II begins. Denote by $x^n$ the real coding the infinite sequence $(x^n_0, x^n_1,\hdots)$.
\item During turn $\omega\cdot n$, player I may decide to advance the game to the second stage, if she has not done so before.
\item If so, then Player I must play some $y\in\mathbb{R}$.
\item If so, Player II must respond with some $m \in\mathbb{N}$.
\end{enumerate} 
Player I wins the game if, and only if, there is some $n$ such that the game advances to the second stage on turn $\omega\cdot n$, and if $y\in\mathbb{R}$ and $m\in\mathbb{N}$ are the two numbers played immediately afterwards, then
$$R(\langle x^0, x^1,\hdots, x^n\rangle\upharpoonright m,x\upharpoonright m,\vec a\upharpoonright m, y\upharpoonright m),$$
where $x^i = (x^i_0, x^i_1,\hdots)$.

\begin{claim}
Player I has a winning strategy in $G(x)$ if, and only if, she has one in $H(x)$.
\end{claim}
\proof
Suppose Player I has a winning strategy $\Sigma$ in $H(x)$. Then she has one in $G(x)$: during a turn $k$ of one of the forms $4n$ or $4n+3$, Player I plays the restriction of $\Sigma$ to the next $\omega$-many moves of the game. Player II's move during turn $4n+1$ or $4n+2$ is a strategy which, when applied to $\Sigma$, yields a real which we will denote by $x^{2n}$ or $x^{2n+1}$, according to the parity of $k$. Since $\Sigma$ is winning for Player I, whenever she plays by the strategy just described, there will be an $n\in\mathbb{N}$ such that Player I can play some $y\in\mathbb{R}$ which ensures that
$$\forall m\, R(\langle x^0, x^1,\hdots, x^n\rangle\upharpoonright m,x\upharpoonright m,\vec a\upharpoonright m, y\upharpoonright m),$$
so the strategy described ensures winning $G(x)$.

Conversely, if Player I has a winning strategy $\Sigma$ in $G(x)$, then she has one in $H(x)$, obtained from using, rather than playing, the strategies provided by $\Sigma$. Since $\Sigma$ is a winning strategy, there can be no infinite play $\langle x^0, x^i, \hdots\rangle$ such that
$\phi(\langle x_0,\hdots, x_n\rangle, x, \vec a)$ holds for no $n\in\mathbb{N}$. Thus, for any play obtained this way, there will be a least $n$ for which $\phi(\langle x_0,\hdots, x_n\rangle, x, \vec a)$ holds. Player I can then decide to advance the game to the second stage during turn $\omega\cdot n$, after which she can play a real number witnessing $\phi$. This proves the claim.
\endproof
Since $H(x)$ is an open game of length $\omega^2$ with moves in $\mathbb{R}$, uniformly for all $x\in\mathbb{R}$, and $x\in A$ precisely when Player I has a winning strategy in $H(x)$, we have shown that $A \in \Game^\mathbb{R}_{\omega^2}\SIGMA^0_1$.
This proves Lemma \ref{LemmaOpenCubed1}.
\endproof

The proof of Lemma \ref{LemmaOpenCubed3} requires some preliminary observations. The first one we isolate as Lemma \ref{LemmaOpenCubed2} below. Lemmata \ref{LemmaOpenCubed3} and \ref{LemmaOpenCubed2} are both easy consequences of the proofs of known theorems, but not of the statements. 
Recall that a \emph{scale} on a set $A\subset\mathbb{R}$ is a sequence $\{\varphi_i: i\in\mathbb{N}\}$ of norms on $A$ with the following property: suppose that $x_0, x_1,\hdots$ are elements of $A$ converging to some $x\in\mathbb{R}$ (in the usual topology) and suppose that for each $i\in\mathbb{N}$, 
\[\lim_{n\in\mathbb{N}} \varphi_i(x_n) = \lambda_i\]
(in the discrete topology). Then,
\begin{enumerate}
\item $x\in A$, and
\item for every $i\in\mathbb{N}$, $\varphi_i(x) \leq \lambda_i$.
\end{enumerate}
A pointclass $\Gamma$ has the \emph{scale property} if every set in $\Gamma$ has a scale $\vec\varphi$ whose norms are uniformly in $\Gamma$, in which case we say that $\vec\varphi$ itself is in $\Gamma$.
We refer the reader to Moschovakis \cite{Mo09} for further background.
\begin{lemma} \label{LemmaOpenCubed2}
Suppose games of length $\omega^2$ on $\mathbb{N}$ with $\Game^\mathbb{R}$-hyperprojective payoff are determined. Then $\ind(\Game^\mathbb{R})$ has the scale property.
\end{lemma}
\proof
The proof is somewhat standard, so we only sketch it.
The key result which we need is one of Martin \cite[Corollary 7.2]{Ma08a} by which if a pointclass $\Gamma$ has the scale property and games of length $\omega^2$ on $\Gamma$ are determined, then $\Game^\mathbb{R}\Gamma$ has the scale property.

Afterwards, the proof is much like the argument of Moschovakis \cite[Exercise 7C.19]{Mo09}\footnote{The argument in Moschovakis \cite{Mo08} based on the next-quantifier representation of inductive sets adapts similarly.} for showing that inductive sets admit inductive scales: one verifies that for every operator $\phi(X)$ in which $X$ appears positively and which is definable in the language $L(\Game^\mathbb{R})$, every stage $\phi^\lambda$ of the inductive definition on $\phi$ admits a scale which is $\Game^\mathbb{R}$-hyperprojective. This is done by induction: the successor step is immediate from Martin's theorem. For the limit step, one assumes that scales $\{\psi^\gamma_i:i\in\mathbb{N}\}$ on $\phi^\gamma$ have been constructed for $\gamma<\lambda$ and ``joins'' them together into a scale $\{\psi^{<\lambda}_i:i\in\mathbb{N}\}$ on $\phi^{<\lambda}$, where $\psi^{<\lambda}_0(x)$ is the least $\gamma<\lambda$ such that $x \in \phi^\gamma$ and 
\[\psi^{<\lambda}_{i+1}(x) = \psi^{\psi^\lambda_0(x)}_i(x)\]
for all $i \in\mathbb{N}$.

This defines a scale on $\phi^{<\lambda}$. However, in order to continue the induction, one must check that this scale is $\Game^\mathbb{R}$-hyperprojective.
This is done as in \cite{Mo09}, using the recursion theorem and the fact that the theorem on the transfer of scales from a set $A$ to $\Game^\mathbb{R}A$ has an effective proof (this is the case because the definitions of the scales are written down in \cite{Ma08a}). 
Finally the same construction, together with the Stage Comparison Theorem, yields a scale on $\phi^\infty$ which is in $\ind(\Game^\mathbb{R})$. This completes the proof of Lemma \ref{LemmaOpenCubed2}.
\endproof

Since every $\Game^\mathbb{R}$-hyperprojective set is given by a $\Game^\mathbb{R}$-inductive definition that takes ${<}|\ind(\Game^\mathbb{R})|$-many steps (by the Boundedness Theorem), the preceding argument yields: 
\begin{corollary}[to the proof] \label{CorOpenCubed2}
Suppose games of length $\omega^2$ on $\mathbb{N}$ with $\Game^\mathbb{R}$-hyperprojective payoff are determined. Then every $\Game^\mathbb{R}$-hyperprojective set has a $\Game^\mathbb{R}$-hyperprojective scale.
\end{corollary}

By our convention for dealing with games of length $\omega$ on $\mathbb{R}$, a run of one these games is identified with the (single) real coding it. Below, we will also require the additional assumption that if $x$ is a real coding an $n$-tuple $(x_0,\hdots, x_{n-1})$ and $k < n$, then the first $k$ digits of $x$ depend only on $x_0, \hdots, x_{k-1}$; similarly for codes of infinite sequences. This has the consequence that, if $(x_0,x_1,\hdots)$ is an infinite sequence, then, letting $y_k$ code $(x_0,\hdots, x_{k-1})$ for all $k\in\mathbb{N}$ and $y$ code $(x_0,x_1,\hdots)$, we have
\[\lim_{k\to\infty} y_k = y,\]
which puts us in a situation in which scales might be of use.

\proof[Proof of Lemma \ref{LemmaOpenCubed3}]
One way of proving this is to adapt the argument for the Third Periodicity Theorem of Moschovakis \cite[Theorem 6E.1]{Mo09} to games on $\mathbb{R}$, though the adaptation is straightforward. We mention the main points. Let $A$ be a $\Game^\mathbb{R}$-hyperprojective set and identify it with the game on $\mathbb{R}$ it induces. Assume without loss of generality that it is Player I who has a winning strategy. By Corollary \ref{CorOpenCubed2}, $A$ has a $\Game^\mathbb{R}$-hyperprojective scale $\{\varphi_i:i\in\mathbb{N}\}$ which by \cite[Lemma 4E.2]{Mo09} we may assume is very good.

Following Moschovakis' proof, define, for even $k$, $W_k$ to be the collection of all $(k+1)$-tuples of reals from which Player I has a winning strategy in $A$ and let $H_k(u,v)$ be the game from \cite[Diagram 6E.1]{Mo09}:
\begin{align*}
\begin{array}{c|cccccc}
I &  &   &     &    a_{k+2}&  & \dots\\
II &  a_{k+1}   &      &     &    & a_{k+3}  &\dots\\\hline
I &      &  b_{k+1}    &   &   &      & \dots\\
II &     &          &  b_{k+2}  &&     &\dots
\end{array}
\end{align*}
Here, $u$ is assumed to be the real $\langle a_0,a_1,\hdots,a_k\rangle$ coding the finite sequence of reals $(a_0,a_1,\hdots,a_k)$ and similarly $v = \langle b_0,b_1,\hdots,b_{k}\rangle$. Each move $a_i$ and $b_i$ is a real number. Let
\begin{align*}
a &= \langle a_0, a_1, a_2, \hdots\rangle, \text{ and}\\
b &= \langle b_0, b_1, b_2, \hdots\rangle,
\end{align*}
Player I wins if, and only if $a \leq_{\varphi_k} b$ and Player II wins if, and only if, $a \not \in A$ or $b <_{\varphi_k} a$.

Let $\psi_k$ be the norm defined as in \cite[Theorem 6E.1]{Mo09}: by setting $u \leq_{\psi_k} v$ if, and only if, Player I has a winning strategy in the game $H_k(u,v)$. Since $\varphi_k$ was a norm from a $\Game^\mathbb{R}$-hyperprojective scale and the $\Game^\mathbb{R}$-hyperprojective sets are closed under the real game quantifier, $\psi_k$ is $\Game^\mathbb{R}$-hyperprojective (uniformly in $k$).

Still following \cite[Theorem 6E.1]{Mo09}, we call a sequence of reals of odd length $(a_0,\hdots,a_k) \in W_k$ \emph{minimal} if 
$$(a_0,\hdots,a_{k-1},a_k) \leq_{\psi_k} (a_0,\hdots,a_{k-1}, b)$$
for all $b\in\mathbb{R}$. 
The point of the norms $\psi_k$ is that it is shown as part of the proof of \cite[Theorem 6E.1]{Mo09} that if $a \in\mathbb{R}^\mathbb{N}$ is an infinite play such that $a \upharpoonright 2k$ is minimal for all $k$, then $a$ is a win for Player I in $A$ (this uses our remark on the convergence of codes before the beginning of this proof). Since $\psi_k$ is $\Game^\mathbb{R}$-hyperprojective, the collection of all minimal (codes for finite) tuples of reals is $\Game^\mathbb{R}$-hyperprojective.

In \cite{Mo09}, where games on $\mathbb{N}$ are being played, the definable winning strategy is obtained by playing the least number which is minimal during each turn. When playing games on $\mathbb{R}$, one cannot do that, but the problem can be fixed in the obvious way:

Since $\Game^\mathbb{R}$-hyperprojective sets are (of course) closed under $\forall^\mathbb{R}$ and have $\Game^\mathbb{R}$-hyperprojective scales, a well-known theorem of Moschovakis \cite{Mo71} implies that every $\Game^\mathbb{R}$-hyperprojective set can be uniformized by a $\Game^\mathbb{R}$-hyperprojective function. Hence, there is a $\Game^\mathbb{R}$-hyperprojective function $f$ such that if $u$ is a tuple of reals of even length in $W_k$, then $f(u) \in\mathbb{R}$ is minimal. The winning strategy for Player I is to play $f(u)$ whenever it is defined.
This completes the proof of the lemma.
\endproof

\section{Further Remarks}
\label{Sectweaker}
By Theorem \ref{mainopencubed}, if open games of length $\omega^3$ are determined, then there is a transitive model of $\kp + \ad_\mathbb{R}$. 
Let us put the hypothesis on open games into context by deriving it from the existence of large cardinals:

\begin{proposition}
Suppose that there are $\omega^2$ Woodin cardinals and a measurable cardinal above them. Then, open games of length $\omega^3$ on $\mathbb{N}$ are determined.
\end{proposition}
\proof
Suppose that there are $\omega^2$ Woodin cardinals and a measurable cardinal above them. It follows from a theorem of Steel \cite{St93} that the extender model $M_{\omega^2}^\sharp$ exists and is  iterable. By work of Neeman \cite[Section 2]{Ne04}, this implies that analytic games of length $\omega^3$ on $\mathbb{N}$ are determined.
\endproof

\begin{corollary}\label{CorMeas}
Suppose that there are $\omega^2$ Woodin cardinals and a measurable cardinal above them. Then, there is a transitive model of $\kp + \ad_\mathbb{R}$ containing $\mathbb{R}$.
\end{corollary}

One of the reviewers has pointed out that the existence of a transitive model of $\kp + \ad_\mathbb{R}$ follows from the existence of $\omega^2$ Woodin cardinals, without the measurable cardinal above. 
We sketch a proof of this.
We need to recall the notion of quasi-determinacy:

\begin{definition}
Let $A \subset \mathbb{R}^\mathbb{N}$. A \emph{winning quasi-strategy} for Player I in $A$ is a function 
\[\Sigma:\bigcup_{n\in\mathbb{N}}\mathbb{R}^{n} \to \mathcal{P}(\mathbb{R})\setminus \{\varnothing\}\]
such that $x \in A$ whenever $x(2n) \in \Sigma(x(0),x(1),\hdots, x(2n-1))$ for all $n\in\mathbb{N}$. A winning quasi-strategy for Player II is defined analogously. We say the game on reals on $A$ is \emph{quasi-determined} if one of the players has a winning quasi-strategy for it.
\end{definition}
A quasi-strategy is essentially a multi-valued strategy. Under some circumstances, quasi-determinacy implies determinacy:
\begin{lemma}\label{LemmaQD}
Let $\Gamma$ and $\Lambda$ be pointclasses and suppose that the following hold:
\begin{enumerate}
\item $\Gamma\subset\Lambda$;
\item $\Gamma$ is closed under continuous preimages and $\Lambda$ is closed under continuous preimages, $\exists^\mathbb{R}$, complements, and finite unions;
\item every game on reals in $\Gamma$ has a winning quasi-strategy in $\Lambda$, and every game on reals in $\Lambda$ is quasi-determined. 
\end{enumerate}
Then, every game on reals in $\Gamma$ is determined. Moreover, if $\Gamma = \Lambda$, then the game has a winning strategy in $\Gamma$.
\end{lemma}
\proof
The lemma is proved by the argument in Larson \cite[Remark 6.1.3]{La17}. We repeat it for convenience, noting the definability constraints. First, observe that every set in $\Lambda$ can be uniformized. To see this, given a set $A \in \Lambda$, consider the game where Player I selects some $x \in\mathbb{R}$ such that $(x,y)\in A$ for some $y$, and Player II plays the digits of such a $y$. Since $\Lambda$ is closed under $\exists^\mathbb{R}$, continuous preimages, complements, and finite unions, this game belongs to $\Lambda$. By quasi-determinacy, Player II has a winning quasi-strategy, from which one can extract a uniformizing function for $A$ (which maps an $x$ as above to the real $y$ whose digits are the lexicographically least allowed by the quasi-strategy).

Afterwards, consider a game on reals with payoff set $A \in \Gamma$. Without loss of generality, suppose that Player I has a winning quasi-strategy for $A$ in $\Lambda$. The set of partial plays of this game consistent with this quasi-strategy belongs to $\Lambda$, so uniformization yields a winning strategy.

Suppose additionally that $\Gamma = \Lambda$. Then, the argument in the first paragraph of the proof shows that the uniformizing function belongs to $\Gamma$, so that the strategy obtained in the preceding paragraph also belongs to $\Gamma$.
\endproof
 
We now sketch the argument for producing a model of $\kp + \ad_\mathbb{R}$ from $\omega^2$ Woodin cardinals.
Suppose there are $\omega^2$ Woodin cardinals and let $M$ be the derived model at their limit (so $M$ is an inner model of a generic extension of $V$). By a theorem of Woodin (see also Trang \cite[Theorem 1.1]{Tr15}) $M$ satisfies $\ad^+$ + $\Theta = \Theta_0$ + ``there is a fine, normal ultrafilter on $\mathcal{P}_{\omega_1}(\mathbb{R})$.'' Work in $M$ and let $\mu$ be the normal, fine ultrafilter on $\mathcal{P}_{\omega_1}(\mathbb{R})$. Suppose $A \subset\mathbb{R}$ is Suslin, co-Suslin in $M$. By a theorem of Martin (see Martin-Steel \cite[Theorem 1.1]{MaSt08a}), $A$ is homogeneously Suslin. Fix homogeneous trees $S$ and $T$ on $\omega\times \kappa$, for some $\kappa<\delta^2_1$, projecting to $A$ and its complement.
Given $\sigma \in\mathcal{P}_{\omega_1}(\mathbb{R})$, Larson \cite[Lemma 12.3.1]{La17} shows that the game on reals on the projection of $S$ is quasi-determined in $L[S,T](\sigma)$. 

Since $S^\sharp$ and $T^\sharp$ exist, we have
\[\sigma = \mathbb{R} \cap L[T,S](\sigma)\] 
for $\mu$-almost every $\sigma \in\mathcal{P}_{\omega_1}(\mathbb{R})$.
For such a $\sigma$, every set in $L[S,T](\sigma)$ is definable in $L[S,T](\sigma)$ from an element of $\sigma$, $S$, $T$, and finitely many ordinal parameters. By fineness and normality, there is a single $x_0 \in\mathbb{R}$ such that for $\mu$-almost all $\sigma$, there is a winning quasi-strategy in $L[S,T](\sigma)$ for the game on reals on the projection of $S$, as computed in $L[S,T](\sigma)$, and such a strategy is definable in $L[S,T](\sigma)$ from $x_0$, $S$, $T$, and ordinal parameters. Since $\hod_{S,T,x_0}^{L[S,T](\sigma)}$ is wellordered, there is a least such strategy, $\Sigma_\sigma$.
A winning quasi-strategy $\Sigma$ for $A$ in $M$ can be obtained by letting $y \in \Sigma(x)$ if, and only if, $y \in \Sigma_\sigma(x)$ for $\mu$-almost all $\sigma$.

Hence, $A$ is quasi-determined in $M$. Moreover, the quasi-strategy for $A$  can be obtained directly from $\mu$, $S$ and $T$, the latter two of which are trees on $\omega\times\kappa$ for some $\kappa<\delta^2_1$, so there is a winning quasi-strategy for $A$ which is Suslin and co-Suslin. Applying Lemma \ref{LemmaQD} with $\Gamma = \Lambda = \DELTA^2_1$, every game on reals with Suslin, co-Suslin payoff in $M$ is determined as witnessed by a Suslin, co-Suslin winning strategy.
Therefore, the companion of the Spector class $(\Sigma^2_1)^M$ satisfies $\kp + \ad_\mathbb{R}$.

With little additional work, we can get (from $\omega^2$ Woodin cardinals alone) a transitive model of $\kp + \ad_\mathbb{R}$ containing $\mathbb{R}$. 
By reflecting to some $V_\lambda$ satisfying a large-enough fragment of $\zfc$, say $\zfc_k$ (i.e., the first $k$ axioms of $\zfc$), and again resorting to Steel \cite{St93}, we see that there is a countable, iterable model $N$ of $\zfc_k$ with $\omega^2$ Woodin cardinals. Letting $M$ denote the derived model of $N$ at the supremum of its Woodin cardinals, the argument above shows that  the companion of $(\Sigma^2_1)^M$ satisfies $\kp + \ad_\mathbb{R}$. Standard arguments show that if $g\subset \coll(\omega,\mathbb{R})$ is generic, then, in $V[g]$, $N$ can be iterated to a model $N'$ whose derived model $M'$ has exactly $\mathbb{R}^V$ as reals. Thus, in $V[g]$, there is a transitive model of $\kp + \ad_\mathbb{R}$ containing $\mathbb{R}^V$ (namely, the companion of $(\Sigma^2_1)^{M'}$). 

A different reviewer has given a short argument for getting such a model in $V$. Assuming without loss of generality that $M$ above is minimal, it is of the form $L_\gamma(M|\lambda^M)$ for some $\gamma$, where $\lambda^M$ is the supremum of the Woodin cardinals of $M$. The derived model of $M$ is then of the form $L_{\gamma'}(\mathbb{R},\mathcal{C})$ for some $\gamma'$, where $\mathcal{C}$ is the closed, unbounded filter on $\mathcal{P}_{\omega_1}(\mathbb{R})$ (see Rodr\'iguez and Trang \cite{RoTr18}); thus, it belongs to $V$.

By varying $k$ (of $\zfc_k$)  above  (provided it is big enough), one obtains different Spector classes of the form $(\Sigma^2_1)^{M'}$ and different models of $\kp+\ad_\mathbb{R}$. The following proposition is proved by relating these models to the one produced in Section \ref{SectPrel}:

\begin{proposition}
Suppose that there are $\omega^2$ Woodin cardinals. Then, open games of length $\omega^3$ are determined.
\end{proposition}
\proof[Proof Sketch]
Choose $k_0$ such that $\zfc_{k_0}$ suffices for the argument above, so that letting $M_{k_0}$ be the associated derived model, $M_{k_0}$ is transitive, contains $\mathbb{R}$, satisfies $\zf_{k_0} + \ad^+$ + ``there is a fine, normal ultrafilter on $\mathcal{P}_{\omega_1}(\mathbb{R})$,'' and the companion of $(\Sigma^2_1)^{M_{k_0}}$ satisfies $\kp + \ad_\mathbb{R}$. Let $k_0<k_1$ and define $M_{k_1}$ analogously. Without loss of generality, there is no inner model of $M_{k_0}$ containing all reals which satisfies $\zf_{k_1} + \ad^+$ + ``there is a fine, normal ultrafilter on $\mathcal{P}_{\omega_1}(\mathbb{R})$.''

The point is that if $M_{k_0}$ and $M_{k_1}$ are as above, then $(\Sigma^2_1)^{M_{k_0}}$ is strictly smaller than $(\Sigma^2_1)^{M_{k_1}}$, and thus every set of reals in $(\Sigma^2_1)^{M_{k_0}}$ is co-Suslin in $M_{k_1}$. Working in $M_{k_0}$, every game on reals with a Suslin, co-Suslin payoff is determined, as witnessed by a Suslin, co-Suslin strategy. Thus, the Suslin, co-Suslin sets are closed under $\Game^\mathbb{R}$. Since $M_{k_0}$ contains all reals, $M_{k_0}$ is correct about a strategy for a game on reals being a winning strategy, so $M_{k_0}$ computes the quantifier $\Game^\mathbb{R}$ correctly when applied to sets that are Suslin, co-Suslin in $M_{k_0}$. Thus, $M_{k_0}$ can carry out inductive definitions with respect to the quantifier $\Game^\mathbb{R}$ and thus contains all $\ind(\Game^\mathbb{R})$ sets.

Depending on what $k_0$ is, an arbitrary $\ind(\Game^\mathbb{R})$ set might or might not be Suslin and co-Suslin in $M_{k_0}$, but it certainly is so in $M_{k_1}$. Since every game on reals with Suslin, co-Suslin payoff in $M_{k_1}$ is determined, as witnessed by a strategy which is Suslin and co-Suslin in $M_{k_1}$, we conclude that all $\ind(\Game^\mathbb{R})$ are determined, as witnessed by a strategy which is Suslin and co-Suslin in $M_{k_1}$. Hence, we can localize the argument of Blass \cite{Bl75} for showing that $\ad_\mathbb{R}$ implies determinacy of all games of length $\omega^2$, as in \cite[Lemma 3.11]{Ag18c}, to show that all games of length $\omega^2$ on $\mathbb{N}$ with $\ind(\Game^\mathbb{R})$ payoff are determined. As observed in Remark \ref{theremark}, this implies that all open games of length $\omega^3$ are determined.
\endproof
A natural question is that of the precise consistency strength of $\kp+\ad_\mathbb{R}$. The results presented here suggest the following conjecture.
Below, recall our convention that $\ad_\mathbb{R}$ (by definition) implies that $\mathbb{R}$ exists.

\begin{conjecture}
The following are equiconsistent:
\begin{enumerate}
\item $\kp$ $+$ $\ad_\mathbb{R}$;
\item\label{clauseWoodins} $\kp$ $+$ ``there are infinitely many limits of Woodin cardinals, cofinal in the ordinals.''
\end{enumerate}
\end{conjecture}

The sentence in item \eqref{clauseWoodins} above can be more precisely interpreted as ``for every ordinal $\alpha$, there are cardinals $\lambda$ and $\delta$, with $\alpha<\lambda<\delta$, such that $V_{\delta}$ exists and satisfies $\zfc$ + `$\lambda$ is a limit of Woodin cardinals.' ''
We also mention the corresponding conjecture for games on natural numbers (see also Remark \ref{remarkad}).

\begin{conjecture}
The following are equiconsistent:
\begin{enumerate}
\item $\kp$ $+$ $\ad$;
\item $\kp$ $+$ ``there are infinitely many Woodin cardinals, cofinal in the ordinals.''
\end{enumerate}
\end{conjecture}

\bibliographystyle{abbrv}
\bibliography{References}

\end{document}